\RequirePackage{fix-cm}
\documentclass[smallextended]{svjour3}       
\smartqed  
\usepackage{graphicx}
\usepackage{amsmath}
\usepackage{graphicx, color}
\usepackage{amscd}
\usepackage{amsfonts}
\usepackage{amssymb}
\usepackage{mathrsfs}
\usepackage{mathtools}
\usepackage [latin1]{inputenc}
 \usepackage{mathptmx}      
\newcommand{\R}{\mathbb{R}}
\newcommand{\N}{\mathbb{N}}

\journalname{Bollettino dell'Unione Matematica Italiana}
\begin{document}

\title{Approximation by regular functions in
Sobolev\\spaces arising from doubly elliptic problems\\
\thanks{The authors  are members of the {\em Gruppo Nazionale per
l'Analisi Ma\-te\-ma\-ti\-ca, la Probabilit\`a e le loro Applicazioni}
(GNAMPA) of the {\em Istituto Nazionale di Alta Matematica} (INdAM).
The manuscript was realized within the auspices of the INdAM -- GNAMPA Projects
{\em Equazioni alle derivate parziali: problemi e mo\-del\-li} (Prot\_U-UFMBAZ-2020-000761).
The first author was also partly supported  by  the {\em Fondo Ricerca di Base di Ateneo -- Eser\-ci\-zio 2017--2019} of the University of Perugia, named {\em PDEs and Nonlinear Analysis, while
the second author by the} {\em Progetto Equazione delle onde con condizioni acustiche,  finanziato  con  il Fondo  Ricerca  di Base, 2019, della Universit\`a degli Studi di Perugia} and by {\em Progetti Equazioni delle onde con condizioni iperboliche ed acustiche al bordo,  finanziati  con  i Fondi  Ricerca  di Base 2017 and 2018, della Universit\`a degli Studi di Perugia}.}}


\titlerunning{Approximation by regular functions}        

\author{Patrizia Pucci  \and Enzo Vitillaro}


\institute{P. Pucci \at
              Dipartimento di Matematica e Informatica,
Universit\`a degli Studi di Perugia\\
Via Vanvitelli 1, 06123 Perugia, Italy \\
              Tel.: +390755855038\\
              Fax: +390755855024\\
              \email{patrizia.pucci@unipg.it}           
           \and
           E. Vitillaro \at
              Dipartimento di Matematica e Informatica,
Universit\`a degli Studi di Perugia\\
Via Vanvitelli 1, 06123 Perugia, Italy\\
              Tel.: +390755855015\\
              Fax: +390755855024\\
              \email{enzo.vitillaro@unipg.it}
}

\date{Received: date / Accepted: date}

\dedication{To the memory of our dear friend Professor Domenico Candeloro with high feelings of admiration for his notable contributions in Mathematics}

\maketitle

\begin{abstract}
The paper deals with a nontrivial density result for $C^m(\overline{\Omega})$ functions, with $m\in\N\cup\{\infty\}$, in the space
$$W^{k,\ell,p}(\Omega;\Gamma)=
\left\{u\in W^{k,p}(\Omega): u_{|\Gamma}\in W^{\ell,p}(\Gamma)\right\},$$
endowed with the norm of $(u,u_{|\Gamma})$ in $W^{k,p}(\Omega)\times W^{\ell,p}(\Gamma)$, where
$\Omega$ is a bounded open subset of $\R^N$, $N\ge 2$, with  boundary $\Gamma$ of class $C^m$, $k\le \ell\le m$ and $1\le p<\infty$.

Such a result is of interest when dealing with doubly elliptic problems involving two elliptic operators, one in $\Omega$ and
the other on $\Gamma$.

Moreover we shall also consider the case when a Dirichlet homogeneous boundary condition is imposed on a relatively open part of $\Gamma$ and, as a preliminary step, we shall prove an analogous result when either $\Omega=\R^N$ or $\Omega=\R^N_+$ and $\Gamma=\partial\R^N_+$.
\keywords{Density results\and  Sobolev spaces \and Smooth functions \and the Laplace--Beltrami operator}
\subclass{46E35 \and 46.38 \and 46M35}
\end{abstract}

\section{Introduction and main results}
\label{intro}

Density results for smooth functions in Sobolev spaces constitute a cornerstone in the classical theory of these spaces and in their applications in PDEs theory. Actually every textbook dealing with Sobolev spaces or PDEs devotes some attention to this subject, see for example \cite{adams,AdamsFournier,brezis2,LeoniSobolev2,lionstata,lionsmagenes1,Mazja,necasbook,triebel}.
\nocite{barcandeloro}
The paper  deals with Sobolev spaces of integer nonnegative order, which are the most classical ones, but the density subject is  standard also when working with Sobolev spaces of fractional order, also known as Sobolev -- Slobodeckij spaces, and with Besov spaces and Bessel -- potential ones. See \cite{triebel}.

When considering Sobolev spaces in $\R^N$, $N\ge 1$, it is well--known that the space of compactly supported smooth functions $C^\infty_c(\R^N)$ is dense in $W^{k,p}(\R^N)$ for $k\in\N$ and $1\le p<\infty$. Since the proof of this result relies on Friedrichs mollifiers and truncation arguments,
following \cite{brezis}, we shall refer to this result as to Friedrichs' Theorem.
When considering the same spaces in an open subset $\Omega$ of $\R^N$, $N\ge 1$, the same arguments  show that $\left\{u\big|_{\Omega}: u\in C^\infty_c(\R^N)\right\}$ is dense in $W^{k,p}(\Omega)$ for $k\in\N$ and $1\le p<\infty$, provided $\Gamma=\partial\Omega$ is continuous (see \cite[Theorem~ 11.35, p. 330]{LeoniSobolev2}).

Moreover, the celebrated Meyers Serrin Theorem (see \cite{MeyersSerrin}) asserts that $C^\infty(\Omega)\cap W^{k,p}(\Omega)$ is dense in $W^{k,p}(\Omega)$, for $k\in\N$ and $1\le p<\infty$, without any regularity assumption on $\Gamma$.
Also this result is treated in most textbooks in Sobolev spaces, see for example \cite{LeoniSobolev2} and \cite{ziemer}.

In the present paper we shall deal with  a Friedrichs' type result. Indeed, as we are going to explain in the sequel, boundary regularity is needed even to state our first main result.

It is mathematical folklore that all density results for smooth functions in Sobolev spaces and in their closed subspaces can be trivially derived by the classical results stated above. On the other hand, for non--closed subspaces, the situation may be different. We refer here to \cite{quarteroni},  where the authors deal with the doubly elliptic problem
$$
\begin{cases} -\Delta u=f \qquad &\text{in
$\Omega$,}\\
\partial_\nu u+\alpha u-\beta\Delta_\Gamma u=h &\text{on $\Gamma$,}
\end{cases}$$
where $\Omega\subset \R^N$, $N\ge 2$,  is a bounded domain with $C^m$ boundary $\Gamma$, $m\in\N\cup\{\infty\}$, $\alpha,\beta$ are positive constants, $\nu$ denotes the outward normal to $\Omega$ and $\Delta_\Gamma$ stands for the Laplace--Beltrami operator on $\Gamma$. In \cite[Remark 2.6]{quarteroni} the authors realized
that the density of $\left\{u\big|_{\Omega}: u\in C^m(\R^N)\right\}$ in the space
\begin{equation}\label{2}
H^m(\Omega;\Gamma)=\left\{u\in H^m(\Omega): u_{|\Gamma}\in H^m(\Gamma)\right\},
\end{equation}
endowed with standard product norm of the couple $(u,u_{|\Gamma})$ in $H^m(\Omega)\times H^m(\Gamma)$, is nontrivial at all. Here and in the sequel $u\mapsto u_{|\Gamma}$ denotes the Trace Operator. See also \cite[Lecture 12]{lionstata}.

The same type of remark was made by the second author of the present paper in \cite{Dresda1,Dresda2}  when dealing with a nonlinear perturbation of the problem
$$
\begin{cases} u_{tt}-\Delta u=0 \qquad &\text{in
$(0,\infty)\times\Omega$,}\\
u=0 &\text{on $(0,\infty)\times\Gamma_0$,}\\
u_{tt}+\partial_\nu u-\Delta_\Gamma u=0 &\text{on $(0,\infty)\times\Gamma_1$,}\\
u(0,x)=u_0(x),\quad u_t(0,x)=u_1(x) &\text{in $\overline{\Omega}$,}
\end{cases}
$$
where $\Gamma$ is of class $C^1$, $\Gamma=\Gamma_0\cup\Gamma_1$ with $\overline{\Gamma_0}\cap\overline{\Gamma_1}=\emptyset$ and $\Gamma_1\not=\emptyset$.
In particular the density of $\left\{u\big|_{\Omega}: u\in C^1(\R^N), \,\, u=0\,\,\text{on}\,\Gamma_0\right\}$
 in the space
\begin{equation}\label{4}
 H^1_{\Gamma_0}(\Omega;\Gamma)=\left\{u\in H^1(\Omega): u_{|\Gamma}\in H^1(\Gamma), u_{|\Gamma}=0\quad\text{on}\,\,\Gamma_0\right\},
\end{equation}
endowed with standard product norm of the couple $(u,u_{|\Gamma})$ in $H^1(\Omega)\times H^1(\Gamma)$,
was remarked to be nontrivial.

Clearly the two results above are particular cases of a slightly more general one. Before formulating it we remark that the space $H^m(\Gamma)$ appearing in \eqref{2} and in \eqref{4} (when $m=1$) is properly defined, through local charts, only when $\Gamma$ is at least $C^{m-1,1}$ (see \cite{grisvard}).

In this paper we shall take $\Gamma$ of class $C^m$ for the sake of simplicity and thus we shall consider $\Omega$ satisfying the following assumption:
\renewcommand{\labelenumi}{{$(H)$}}
\begin{enumerate}
\item {\em $\Omega$ is a bounded open subset of $\R^N$, $N\ge 2$, with $\Gamma=\partial\Omega$ of class $C^m$,
$m\in\N\cup\{\infty\}$. Moreover $\Gamma=\Gamma_0\cup\Gamma_1$, $\overline{\Gamma_0}\cap
\overline{\Gamma_1}=\emptyset$ and $\Gamma_1\not=\emptyset$.}
\end{enumerate}
\noindent
We shall consider, for $k, \ell\in\mathbb N$
 and  $p\in [1,\infty)$,  the Banach space
 $$
W^{k,\ell,p}_{\Gamma_0}(\Omega;\Gamma)=
\left\{u\in W^{k,p}(\Omega): u_{|\Gamma}\in W^{\ell,p}(\Gamma)\,:\,
u_{|\Gamma}=0\mbox{ on }\Gamma_0\right\},
$$
with the standard contraction  $W^{k,\ell,p}_{\Gamma_0}(\Omega;\Gamma)=W^{k,\ell,p}(\Omega;\Gamma)$ when $\Gamma_0$ is empty.

Clearly this space is of some interest only when  $k\le\ell$, since when $\ell<k$ it reduces to $W^{k,p}(\Omega)$ by the Trace Theorem. In the sequel we shall take $k\le\ell\le m$ and  we shall identify $W^{k,\ell,p}_{\Gamma_0}(\Omega;\Gamma)$ with its isometric copy
$$
W^{k,\ell,p}_{\Gamma_0}(\Omega;\Gamma)=
\left\{(u,v)\in W^{k,p}(\Omega)\times W^{\ell,p}(\Gamma)\,:\,
v=u_{|\Gamma},\,\,v=0\mbox{ on }\Gamma_0\right\}.
$$
We shall endow it with the norm inherited from the product space. Let us remark that the case $\Gamma_0=\emptyset$ is also included in our treatment.
The first main result is

\begin{theorem}\label{theorem1}
If assumption $(H)$ holds then
$$Y=\left\{u\big|_{\Omega}\,:\,u\in C^m_c(\mathbb R^N),\,\,u=0\mbox{ on }\Gamma_0\right\}$$
is dense in $W^{k,\ell,p}_{\Gamma_0}(\Omega;\Gamma)$ for all
$k$, $\ell\in\mathbb N$, with  $k\le\ell\le m$, and any $p\in[1,\infty)$.
\end{theorem}
 The proof of Theorem~\ref{theorem1} relies on the combination of the standard localization technique with an analogous result in the case $\Omega=\R^N$, $\Gamma_0=\emptyset$ and $\Gamma_1=\partial\R^N_+$.
 Since this result could be of some independent interest we shall state here as our second main result.
 To state it we introduce the standard notation $x=(x',x_N)\in \mathbb{R^N}$, with $x'\in\mathbb R^{N-1}$ and~$x_N\in\mathbb R$,
 \begin{gather*}
 B_r(\mathbb R^{N-1})=\{x'\in\mathbb R^{N-1}\,:\,|x'|<r\}\quad\mbox{for any $r>0$,}\\
 \mathbb R^N_+=\{x=(x',x_N)\in\mathbb R^N\,:\,x_N>0\},\quad
\partial\mathbb R^N_+=\{x=(x',0)\in\mathbb R^N\,:\,x'\in\mathbb R^{N-1}\},
\end{gather*}
and  for $k$, $\ell\in\mathbb N$ and $p\in[1,\infty)$,  the Banach space
\begin{equation}\label{7}
W^{k,\ell,p}(\R^N;\partial\mathbb R^N_+)=
\left\{u\in W^{k,p}(\R^N): u_{|\partial\R^N_+}\in W^{\ell,p}(\partial\R^N_+)\right\},
 \end{equation}
where $W^{\ell,p}(\partial\mathbb R^N_+)$
is naturally identified with $W^{\ell,p}(\mathbb R^{N-1})$. According to the previous identification, we shall also identify it with
$$W^{k,\ell,p}(\mathbb R^N;\partial\mathbb R^N_+)=
\{(u,v)\in W^{k,p}(\mathbb R^N)\times W^{\ell,p}(\partial\mathbb R^N_+)\,:\,
v=u_{|\partial\mathbb R^N_+}\},$$
and we shall endow it with the norm inherited from the product space. Also in this case only the case $k\le \ell$ is of some interest.
The second main result is

\begin{theorem}\label{theorem2}
Let $k$, $\ell\in\mathbb N$ and $p\in[1,\infty)$. For any
$u\in W^{k,\ell,p}(\mathbb R^N;\partial\mathbb R^N_+)$
there exists a sequence $(u_n)_n$ in $W^{k,\ell,p}(\mathbb R^N;\partial\mathbb R^N_+)\cap
C^\infty(\mathbb R^N)$ such that $\mathrm{supp}\,u_n\subseteq\mathrm{supp}\,u
+\overline{B_{1/n}(\mathbb R^N)}$ for all $n$ and $u_n\to u$ in~$W^{k,\ell,p}(\mathbb R^N;\partial\mathbb R^N_+)$.
\end{theorem}
Since Theorem~\ref{theorem2} does not  look as the exact translation of Theorem~\ref{theorem1} in the case $\Omega=\R^N$, we would like to remark a trivial consequence of it. To state it
we set, for $k$, $\ell\in\mathbb N$ and $p\in[1,\infty)$,  the Banach space
\begin{equation}\label{8}
W^{k,\ell,p}(\R^N_+;\partial\mathbb R^N_+)=
\left\{u\in W^{k,p}(\R^N_+): u_{|\partial\R^N_+}\in W^{\ell,p}(\partial\R^N_+)\right\},
 \end{equation}
 endowed with the norm of the couple $\left(u, u_{|\partial\R^N_+}\right)$ in the product space.
 Since any element of $W^{k,p}(\R^N_+)$ possesses an extension in $W^{k,p}(\R^N)$ (see for example \cite[Theorem 5.19 p.148]{adams}), by Theorem~\ref{theorem2} we immediately get the next result.
 \begin{corollary} Let $k$, $\ell\in\mathbb N$ and $p\in[1,\infty)$. Then
 $$\left\{ u\big|_{\R^N_+} : u\in C^\infty(\R^N)\right\}\cap W^{k,\ell,p}(\R^N_+;\partial\mathbb R^N_+)$$
 is dense in $W^{k,\ell,p}(\R^N_+;\partial\mathbb R^N_+)$.
 \end{corollary}
  The proof of Theorem~\ref{theorem2} is based on  identifying the space $W^{k,p}(\R^N)$ with its vectorial version
 \begin{equation}\label{Xkp}
 X^{k,p}=\bigcap_{j=0}^k W^{j,p}(\mathbb R;W^{k-j}(\mathbb R^{N-1})),
 \end{equation}
 on choosing mollifiers in separate form in the variables
 $x'\in\mathbb R^{N-1}$ and $x_N\in\mathbb R$ and appropriately selecting their support radii.

 Since a similar identification does not look to be trivial for non integer values of $k$, the extension of Theorems~\ref{theorem1} and \ref{theorem2} to non integer values of $k$ and $\ell$ is not immediate.
  In the next section we are going to give the  proofs of Theorems~\ref{theorem1}--\ref{theorem2}.
\section{Proofs}

\noindent
{\em Proof of Theorem~$\ref{theorem2}$.}
Let $(\rho'_n)_n$ and $(\rho''_n)_n$ be two sequences of standard mollifiers
in~$\mathbb R^{N-1}$ and in~$\mathbb R$, respectively. That is
\begin{alignat*}5
&\rho'_n\in C^\infty_c(\mathbb R^{N-1}),\quad
&&\mathrm{supp}\,&&\rho'_n\subseteq\overline{B_{1/n}(\mathbb R^{N-1})},
\quad&&\int_{\mathbb R^{N-1}}\rho'_ndx'=1,\quad &&\rho'_n\ge0
\mbox{ in }\mathbb R^{N-1},\\
&\rho''_n\in C^\infty_c(\mathbb R),\quad
&&\mathrm{supp}\,&&\rho''_n\subseteq[-1/n,1/n],
\quad&&\int_{\mathbb R}\rho''_ndx_N=1,\quad &&\rho''_n\ge0
\mbox{ in }\mathbb R.
\end{alignat*}
Let $(\rho_{m,n})_{m,n}$ be the double sequence of smooth
functions in~$\mathbb R^N$ defined in the separate form
$$\rho_{m,n}(x)=\rho'_m(x')\,\rho''_n(x_N)\quad\mbox{for }
x=(x',x_N)\in\mathbb R^N.$$
Thus, for any couple of strictly increasing sequences
$(\sigma_n)_n$, $(\tau_n)_n$ in $\mathbb N$ such that
$\sigma_n$, $\tau_n\ge2n$ for all $n$ the sequence
$(\rho_n)_n$, with $\rho_n=\rho_{\sigma_n,\tau_n}$,
is a standard mollifying sequence in~$\mathbb R^N$.

Fix $u\in W^{k,\ell,p}(\mathbb R^N;\partial\mathbb R^N_+)$
and a couple of strictly increasing sequences
$(\sigma_n)_n$, $(\tau_n)_n$ in $\mathbb N$ such that
$\sigma_n$, $\tau_n\ge2n$ for all $n$, which we select lately.
Put $u_{m,n}=\rho_{m,n}\ast u$. Standard properties of convolution
and mollifiers, cf. Propositions~4.18 and~4.20 of~\cite{brezis2}, show that
$\mathrm{supp}\,u_{\sigma_n,\tau_n}\subseteq\mathrm{supp}\,u+
\overline{B_{1/n}(\mathbb R^N)}$
since $\sigma_n$, $\tau_n\ge2n$ for all $n$,
and $u_{\sigma_n,\tau_n}\in C^\infty(\mathbb R^N)$.
Moreover, $u_{\sigma_n,\tau_n}\to u$ in~$W^{k,\ell,p}(\mathbb R^N,\partial\mathbb R^N_+)$
by Lemma~9.1 of~\cite{brezis2}. To complete the proof it is then enough to conveniently choose
$(\sigma_n)_n$, $(\tau_n)_n$ in such a way that
\begin{equation}\label{bound}
u_{\sigma_n,\tau_n}\big|_{\partial\mathbb R^N_+} \to u\big|_{\partial\mathbb R^N_+}
\quad\mbox{in }W^{\ell,p}(\partial\mathbb R^N_+)\quad\text{as $n\to\infty$}.
\end{equation}
To prove \eqref{bound} we first note that
$W^{k,p}(\mathbb R^N)$ is algebraically and topologically
isomorphic to its vectorial version $X^{k,p}$ defined in \eqref{Xkp},
endowed with the standard norm
$$\|\overrightarrow{u}\|_{X^{k,p}}=\left(\sum_{j=0}^k\|u\|_{W^{j,p}(\mathbb R;W^{k-j}(\mathbb R^{N-1}))}^p\right)^{1/p}\quad\mbox{for any }\overrightarrow{u}\in  X^{k,p},$$
via the isomorphism which associates to each $u\in W^{k,p}(\mathbb R^N)$
its vector--valued version $\overrightarrow{u}$ defined by
$$\overrightarrow{u}(x_N)=u(\cdot,x_N)\quad\mbox{for all }x_N\in\mathbb R.$$

This result, which is well--known (see for example \cite[Example~7.34]{adams} or
\cite[Example~3, p.490]{dautraylionsvol5} in the similar case of $\R^N_+$), can be easily proved.
Indeed the case $k=1$ follows by the general theory of vector--valued Sobolev functions (see \cite[Appendix~A]{breziscazenave}) and the extension of the generalized Leibnitz formula \cite[Theorem~3, p. 303]{Evans} to the duality product in Banach spaces and to  $1\le p<\infty$. The general case then follows by induction.

Furthermore, $u\big|_{\partial\mathbb R^N_+}=\overrightarrow{u}(0)$
for any $u\in W^{k,p}(\mathbb R^N)$ thanks to the identification
of~$W^{k,p}(\mathbb R^{N-1})$ with $W^{k,p}(\partial\mathbb R^N_+)$.
Actually, the equality $u\big|_{\partial\mathbb R^N_+}=\overrightarrow{u}(0)$ is true for continuous functions and then by
density in the entire $W^{k,p}(\mathbb R^N)$. Hence
$W^{k,\ell,p}(\mathbb R^N,\partial\mathbb R^N_+)$
is isomorphic by the same identification to its vectorial version
$$X^{k,\ell,p}=\{(\overrightarrow{u},v)\in
X^{k,p}\times W^{\ell,p}(\mathbb R^{N-1})\,:\,v=\overrightarrow{u}(0)\},$$
endowed with the product norm
$$\|\overrightarrow{u}\|_{X^{k,\ell,p}}=\left(\|\overrightarrow{u}\|_{X^{k,p}}^p
+\|\overrightarrow{u}(0)\|_{W^{\ell,p}(\mathbb R^{N-1})}^{p}\right)^{1/p}.$$
Consequently, to show \eqref{bound} is equivalent to prove that
\begin{equation}\label{vect}
\overrightarrow{u}_{\sigma_n,\tau_n}(0) \to \overrightarrow{u}(0)
\quad\mbox{in }W^{\ell,p}(\mathbb R^{N-1}),
\end{equation}
where $(\sigma_n)_n$, $(\tau_n)_n$  will be chosen later.
In order to prove \eqref{vect} let us denote by $\ast'$ and $\ast''$
the convolution in~$\mathbb R^{N-1}$ and in~$\mathbb R$, respectively.

Since $u\in W^{k,\ell,p}(\mathbb R^N,\partial\mathbb R^N_+)$,
then $\overrightarrow{u}\in X^{k,\ell,p}$ and for all $n$
set
$$\overrightarrow{u}_n=\rho_n'\ast'\overrightarrow{u}.$$
Now Proposition 4.20 and Theorem~4.15 of~\cite{brezis2}
imply that for any $\rho'\in C^\infty_c(\mathbb R^{N-1})$,
any $v\in L^p(\mathbb R^{N-1})$ and any multi--index
$\alpha=(\alpha_1,\dots,\alpha_{N-1})$
$$D^\alpha(\rho'\ast'v)=(D^\alpha\rho')\ast'v,$$
so that $\rho'\ast'v\in W^{\ell,p}(\mathbb R^{N-1})$ and
$$\|D^\alpha(\rho'\ast'v)\|_p\le\|D^\alpha\rho'\|_1\|v\|_p,$$
whenever $|\alpha|\le\ell$. Consequently, the linear operator
$v\mapsto \rho'\ast'v$ is bounded from $L^p(\mathbb R^{N-1})$
into $W^{\ell,p}(\mathbb R^{N-1})$.

The continuity of the embeddings $X^{k,p}\hookrightarrow
W^{1,p}(\mathbb R;L^p(\mathbb R^{N-1}))\hookrightarrow
C_b(\mathbb R;L^p(\mathbb R^{N-1}))$ and the fact that
$\overrightarrow{u}\in X^{k,p}$ imply that
$\overrightarrow{u}_n\in C(\mathbb R;W^{\ell,p}(\mathbb R^{N-1}))$.
Now $\overrightarrow{u}(0)=u\big|_{\partial\mathbb R^N_+}$
is in $W^{\ell,p}(\mathbb R^{N-1})$ so that as $n\to\infty$
\begin{equation}\label{I}
\overrightarrow{u}_n(0)\to\overrightarrow{u}(0)\quad
\mbox{in }W^{\ell,p}(\mathbb R^{N-1}).
\end{equation}
Fubini's theorem yields for a.e. $x=(x',x_N)\in\mathbb R^N$
that
\begin{align*}
u_{m,n}(x)&=\int_{\mathbb R^N}\rho_m'(x'-y')\rho_n''(x_N-y_N)
u(y',y_N)dy'dy_N\\
&=\int_{\mathbb R}\rho_n''(x_N-y_N)\big((\rho_m'\ast'\overrightarrow{u})(x')\big)
(y_N)dy_N\\
&=\big(\rho_n''\ast''(\rho_m'\ast'\overrightarrow{u})(x')\big)(x_N).
\end{align*}
Hence $\overrightarrow{u}_{m,n}=\rho_n''\ast''(\rho_m'\ast'\overrightarrow{u})=
\rho_n''\ast''\overrightarrow{u}_n$. Since $\overrightarrow{u}_n\in
C(\mathbb R;W^{\ell,p}(\mathbb R^{N-1})$ the trivial extension to the vectorial case of well--known properties of regularization, that is of \cite[Proposition~4.2]{brezis2}, yield that
for $m$ fixed $\overrightarrow{u}_{m,n}\to\overrightarrow{u}_{m}$
in~$W^{\ell,p}(\mathbb R^{N-1})$ as $n\to\infty$, uniformly on
compact sets of~$\mathbb R$. Consequently,
\begin{equation}\label{II}
\overrightarrow{u}_{m,n}(0)\to\overrightarrow{u}_{m}(0)\quad\mbox{as }
n\to\infty.
\end{equation}
Combining \eqref{I} with \eqref{II} for $\sigma_n=2n$ we get that
$\overrightarrow{u}_{2n}(0)\to\overrightarrow{u}(0)$
in~$W^{\ell,p}(\mathbb R^{N-1})$ and that for any $n$ there exists
$\tau_n\ge 2n$, with $(\tau_n)_n$ strictly increasing, such that
$$\|\overrightarrow{u}_{2n,\tau_n}(0)-\overrightarrow{u}_{2n}(0)\|_
{W^{\ell,p}(\mathbb R^{N-1})}<1/n.$$
Thus, $\overrightarrow{u}_{2n,\tau_n}(0)\to \overrightarrow{u}(0)$
in~$W^{\ell,p}(\mathbb R^{N-1})$, that is \eqref{vect} holds.
Finally this
completes the proof.
\qed
\bigskip

\noindent{\em Proof of Theorem~$\ref{theorem1}$.}
Let us start  by fixing  some usual notation. We set
\begin{gather*}
Q=B_1(\mathbb R^{N-1})\times(-1,1),\quad Q_+=Q\cap\mathbb R^N_+,\quad
Q_0=Q\cap\partial\mathbb R^N_+.
\end{gather*}
Moreover we shall denote by $u_{|\Gamma_i}$ the restriction of $u_{|\Gamma}$ to $\Gamma_i$
 for $i=0,1$.
From the assumption that  $\overline{\Gamma_0}\cap
\overline{\Gamma_1}=\emptyset$ it follows that $\Gamma_0$ and $\Gamma_1$ are
compact. Then, using the definition of $C^m$ regular open set, see Chapter~9
of~\cite{brezis2}, there are open subsets $V_1,\dots,V_r$, $V_{r+1},\dots, V_s$
in~$\mathbb R^N$ such that
$$\Gamma_0\subseteq\bigcup_{j=1}^rV_j,\quad \Gamma_1\subseteq\bigcup_{j=r+1}^sV_j,
\quad V_j\cap\Gamma_1=\emptyset\mbox{ if }j=1,\dots,r,\quad
V_j\cap\Gamma_0=\emptyset\mbox{ if }j=r+1,\dots,s,$$
and bijiective maps $H_j:Q\to V_j$, $j=1,\dots,s$, such that
\begin{gather*}H_j\in C^m(\overline Q),\quad H_j^{-1}\in C^m(\overline{V_j}),\quad
H_j(Q_+)=V_j\cap\Omega,\quad H_j(Q_0)=V_j\cap\Gamma_0,\,\,
j=1,\dots,r,\\
H_j(Q_0)=V_j\cap\Gamma_1,\,\,
j=r+1,\dots,s.
\end{gather*}
Moreover, from Lemma~9.3 of~\cite{brezis2}, see also \cite[Theorem~C.21]{LeoniSobolev2}
for a proof, there are functions $\theta_0,\dots,\theta_s\in
C^\infty(\mathbb R^N)$ such that
$$0\le\theta_j\le1,\quad\sum_{j=0}^s\theta_j=1\mbox{
in }\mathbb R^N,\quad\theta_0\big|_{\Omega}\in C^\infty_c(\Omega),
\quad\mathrm{supp}\,\theta_j\subset\subset V_j\mbox{ for }
j=1,\dots,s.$$
Now let $u\in W^{k,\ell,p}(\Omega,\Gamma)$. Since
$u\in W^{k,p}(\Omega)$ and $\Omega$ is $C^m$ regular
we can extend it to $\tilde u\in W^{k,p}(\mathbb R^N)$
by Theorem~4.26 of~\cite{adams}. Put $u_j=\theta_j\tilde u\in
W^{k,p}(\mathbb R^N)$, so that
\begin{gather*}\tilde u=\sum_{j=0}^s u_j\mbox{ in }\mathbb R^N,\quad
u_{|\Gamma_0}=\sum_{j=0}^r{u_j}_{|\Gamma_0},\quad
u_{|\Gamma_1}=\sum_{j=r+1}^s{u_j}_{|\Gamma_1},\\
\mathrm{supp}\,u_0\subset\subset\Omega,\quad
\mathrm{supp}\,u_j\subset\subset V_j\mbox{ for }
j=1,\dots,s.
\end{gather*}
We shall show that each $u_j$, $j=0,\dots,s$, can be approximated
by elements of $Y$ in the $W^{k,\ell,p}_{\Gamma_0}(\Omega;\Gamma)$
norm. For $j=0$ there is nothing to prove since by Lemma~3.15
of~\cite{adams} there is a sequence $(\psi_n^0)_n$ in~$C^\infty_c(\Omega)$
such that $\psi_n^0\to u_0\big|_\Omega$ in~$W^{k,p}(\Omega)$
and ${\psi_n^0}_{|\Gamma}={u_0}_{|\Gamma}=0$ for all $n$.
Indeed, ${\psi_n^0}_{|\Gamma}=0$ by \cite[Theorem~15.29 p. 475]{LeoniSobolev2} and so
${u_0}_{|\Gamma}=0$ by the continuity of the Trace Operator and the fact
that $\psi_n^0\to u_0\big|_\Omega$ in~$W^{k,p}(\Omega)$.

Fix now $j=1,\dots,s$. By Theorem~3.41 of~\cite{AdamsFournier}
on the stability of Sobolev spaces with respect to coordinate
transformations the linear operator $v\mapsto v\cdot H_j$ is bounded from
$W^{k,p}(\Omega)$ onto~$W^{k,p}(\Omega)$, with bounded inverse.
Hence $u_j\cdot H_j\in W^{k,p}(Q)$. Moreover, $\mathrm{supp}\,(u_j\cdot H_j)\subset\subset Q$,
since $\mathrm{supp}\,u_j\subset\subset V_j$. Then Lemma~3.15 of~\cite{adams}
and Theorem~5.29 of~\cite{AdamsFournier}  yield that $u_j\cdot H_j\in W^{k,p}_0(Q)$ and
so its trivial extension $v_j=\widetilde{u_j\cdot H_j}$ to the whole of $\mathbb R^N$
belongs to~$W^{k,p}(\mathbb R^N)$.

By the definition of Sobolev spaces on $\Gamma$ we have
$\widetilde{u_j\cdot H_j}\big|_{\partial\mathbb R^N_+}\in
W^{\ell,p}(\mathbb R^{N-1})$. Thus
$v_j\in W^{k,\ell,p}(\mathbb R^{N};\partial\mathbb R^N_+)$.
By Theorem~\ref{theorem2} there exists a sequence
$$(\varphi_{j,n})_n\quad\text{
in\quad  $W^{k,\ell,p}(\mathbb R^{N},\partial\mathbb R^N_+)\cap
C^\infty(\mathbb R^N)$}$$
 such that $\varphi_{j,n}\to v_j$
in~$W^{k,\ell,p}(\mathbb R^{N};\partial\mathbb R^N_+)$ as $n\to\infty$ and
$\mathrm{supp}\,\varphi_{j,n}\subseteq\mathrm{supp}\,v_j+
\overline{B_{1/n}(\mathbb R^N)}$. But $\mathrm{supp}\,v_j\subset\subset Q$
so that there exists $r_j\in(0,1)$ such that $\varphi_{j,n}\in
C^\infty_c(Q)$ and $\mathrm{supp}\,\varphi_{j,n}\subseteq
\overline{B_{r_j}(\mathbb R^{N-1})}\times[-r_j,r_j]$ for $n$ sufficiently large.

Let us now distinguish two cases: $j\in\{r+1,\dots,s\}$
and $j\in\{1,\dots,r\}$. When $j\in\{r+1,\dots,s\}$
we set $\psi_{j,n}=\varphi_{j,n}\cdot H_j^{-1}\in C^m_c(V_j)$.
An application of Theorem~3.41 of~\cite{AdamsFournier} yields that
$\psi_{j,n}\to u_j$ in~$W^{k,p}(\mathbb R^{N})$ as $n\to\infty$,
while ${\psi_{j,n}}_{|\Gamma}\to {u_j}_{|\Gamma}$
in~$W^{\ell,p}(\Gamma)$ by the definition of~$W^{\ell,p}(\Gamma)$.
But $\psi_{j,n}\in W^{k,\ell,p}_{\Gamma_0}(\Omega;\Gamma)$, since
in this case $V_j\cap\Gamma_0=\emptyset$. Thus
$\psi_{j,n}\to u_j$ in~$W^{k,\ell,p}_{\Gamma_0}(\Omega;\Gamma)$,
as stated.

When $j\in\{1,\dots,r\}$ we do not know any longer that
$\varphi_{j,n}\cdot H_j^{-1}$ vanishes on $\Gamma_0$, so that
we need to conveniently modify $\varphi_{j,n}$. To this aim
we introduce a cut--off function $\xi\in C^\infty(\mathbb R)$ such that
$\xi(0)=1$ and $\mathrm{supp}\,\xi\subseteq[-r_j,r_j]$. The
linear operator $L:W^{k,p}(\mathbb R^{N-1})\to W^{k,p}(\mathbb R^{N})$
defined by $(Lw)(x)=w(x')\xi(x_N)$ for all $w\in C^\infty(\mathbb R^{N-1})$
is bounded and by density $L$ is well defined in the entire $W^{k,p}(\mathbb R^{N-1})$.
Clearly, $Lw$ is of class $C^\infty(\mathbb R^{N})$ and $Lw\big|_{\partial\mathbb R^N_+}=w$
for all $w\in C^\infty(\mathbb R^{N-1})$, while
$Lw\in C^\infty_c(\mathbb R^{N})$ for all $w\in C^\infty_c(B_1(\mathbb R^{N-1}))$.
Set $\chi_{j,n}=L{\varphi_{j,n}}_{|\partial\mathbb R^N_+}$, so that
${\chi_{j,n}}_{|\partial\mathbb R^N_+}={\varphi_{j,n}}_{|\partial\mathbb R^N_+}$ and moreover
${\chi_{j,n}}_{|\partial\mathbb R^N_+}\to0$ in~$W^{\ell,p}(\mathbb R^{N-1})$
and so in~$W^{k,p}(\mathbb R^{N-1})$, since $k\le \ell$. The fact that $L$ is bounded implies that ${\chi_{j,n}}\to0$ in~$W^{k,p}(\mathbb R^N)$.
Put $\tilde{\varphi}_{j,n}=\varphi_{j,n}- \chi_{j,n}$. Then
$\tilde{\varphi}_{j,n}\in C^\infty_c(\mathbb R^N)$
for all $n$ and $\tilde{\varphi}_{j,n}\to v_j$
in $W^{k,\ell,p}(\mathbb R^{N};\partial\mathbb R^N_+)$,
$\mathrm{supp}\,\tilde{\varphi}_{j,n}\subset\subset Q$ and
finally $\tilde{\varphi}_{j,n}=0$
on~$\partial\mathbb R^N_+$. Set $\psi_{j,n}=\tilde{\varphi}_{j,n}\cdot H_j^{-1}\in
C^m_c(V_j)$. Consequently, as in the previous case, $\psi_{j,n}\to w_j$
in~$W^{k,\ell,p}_{\Gamma_0}(\Omega;\Gamma)$.

Lastly, the previous steps
show that $\psi_n=\sum_{j=0}^{s}\psi_{j,n}$, which is in $Y$,
converges to $u$ in~$W^{k,\ell,p}_{\Gamma_0}(\Omega;\Gamma)$.
This completes the proof.
\qed


\def\cprime{$'$}


\end{document}